\documentclass[letterpaper,10pt,conference]{ieeeconf}
\pdfoutput=1 

\IEEEoverridecommandlockouts
\usepackage{cite,amsmath,amssymb,graphicx,color}
\newcommand{\ones}{\mathbf{1}}

\newcommand{\tr}{\mathrm{Tr}}
\newcommand{\ignore}[1]{}

\title{\LARGE \bf A Majorization-Minimization Approach to\\ Design of Power
Transmission Networks}

\author{Jason K. Johnson and Michael Chertkov%
\thanks{J.~Johnson and M.~Chertkov are both with the Center for Nonlinear
Studies and Theoretical Division T-4 of Los Alamos National Laboratory, Los
Alamos, NM 87544.  M.~Chertkov is also affiliated with the New Mexico
Consortium, Los Alamos, NM 87544. {\tt\small jasonj@lanl.gov,
chertkov@lanl.gov}}%
}

\begin{document}
\maketitle
\thispagestyle{empty}
\pagestyle{empty}

\begin{abstract}
We propose an optimization approach to design cost-effective electrical power
transmission networks.  That is, we aim to select both the network structure and
the line conductances (line sizes) so as to optimize the trade-off between
network efficiency (low power dissipation within the transmission network) and
the cost to build the network. We begin with a convex optimization method based
on the paper ``Minimizing Effective Resistance of a Graph'' [Ghosh, Boyd \&
Saberi]. We show that this (DC) resistive network method can be adapted to the
context of AC power flow.  However, that does not address the combinatorial
aspect of selecting network structure.  We approach this problem as selecting a
subgraph within an over-complete network, posed as minimizing the (convex)
network power dissipation plus a non-convex cost on line conductances that
encourages sparse networks where many line conductances are set to zero.  We
develop a heuristic approach to solve this non-convex optimization problem
using: (1) a continuation method to interpolate from the smooth, convex problem
to the (non-smooth, non-convex) combinatorial problem, (2) the
majorization-minimization algorithm to perform the necessary intermediate smooth
but non-convex optimization steps.  Ultimately, this involves solving a sequence
of convex optimization problems in which we iteratively reweight a linear cost
on line conductances to fit the actual non-convex cost. Several examples are
presented which suggest that the overall method is a good heuristic for network
design.  We also consider how to obtain sparse networks that are still robust
against failures of lines and/or generators.
\end{abstract}

\section{Introduction}

The power grid of today was not systematically planned but grew in a piecemeal fashion. 
In spite of this it is largely reliable, arguably
among the greatest engineering achievements of the 20th century. However, this
status quo is now challenged with increased demand and stress on the aging
network leading to extremely costly and growing-in-scale blackouts and
operational problems. A shift towards renewable sources of energy will further
stress the grid as these resources are intermittent and thus not reliable in
the traditional sense.  These changes emphasize the importance of incorporating new and
extending existing infrastructure in a systematic way. In this paper we present
a proof of principles study suggesting an efficient algorithmic approach for
optimal or close to optimal power grid design.

\subsection{Motivation}

A key challenge in updating and extending the power grid is determining where to
place new transmission, generation and storage facilities or in some cases how
to design a new grid from scratch.  Specifically, the present theoretical study
was motivated by the national challenge of integrating renewables into operation
of the existing US grid. Renewable generation, such as wind and solar, are
intermittent. Moreover, regions where wind is plentiful often lack adequate
transmission lines. Effective and reliable exploitation of renewables requires
planning. The National Renewable Energy Laboratory's (NREL) WinDS project
\cite{20perwind,sbhlong03} is an excellent first step, however, it does not
account for power flow stability or grid resiliency. A study of the WinDS
solution performed at LANL \cite{LorenAlanRussell} has discovered that it
results in an often infeasible electric grid suggesting there is a problem in
generating globally optimal solutions that accommodate intermittent renewable
generation.

Our paper develops an approach towards the challenging problem of planning
cost-effective and robust extensions of the power grid to accommodate growing
demand and long-term addition of renewables.  Our approach may also provide a
starting point for practical planning approaches such as the one proposed in
\cite{LorenAlanRussell}.

\subsection{Related Work}

The initial inspiration for our approach was the convex network optimization
methods of Ghosh, Boyd and Saberi \cite{gbs08}. Building on earlier work
\cite{bvg01}, they consider the problem of minimizing the total resistance of an
electrical network subject to a linear budget on line conductances, where they
interpret the total resistance metric as the expected power dissipation within
the network under a random current model.  We extend their work by also
selecting the network structure. We impose sparsity on that structure in a
manner similar to a number of methods that modify a convex optimization problem
by adding some non-convex regularization to obtain sparser solutions, such as in
compressed sensing \cite{cwb08,chartrand07,chartrand08} or edge-preserving image
restoration \cite{idier08}. The method of Candes \emph{et al} \cite{cwb08} is
especially relevant to our approach.  They recommend the
majorization-minimization algorithm \cite{hl04} as a heuristic approach to
sparsity-favoring non-convex optimization. 

Another important element of our
approach is that we follow a similar strategy as in the graduated non-convexity
algorithm \cite{bz87} in that we solve a sequence of optimization
problems that interpolates from a convex relaxation of the actual non-convex
problem. A somewhat similar approach has been used to obtain sparse transport
networks \cite{bm07}.

\subsection{Our Contributions}

\begin{itemize}

\item We adapt the convex network optimization approach of Ghosh \emph{et al}
\cite{gbs08} to design power transmission networks by demonstrating how AC power
flow (to first order approximation) can be modeled by a (DC) resistive network
model and specializing Ghosh \emph{et al}'s random current model and linear cost
on lines to fit our application.

\item We propose a non-convex, discontinuous generalization of this problem that
more strongly encourages sparsity in the network solution by adding a fixed cost
for each (non-zero conductance) line.  We develop a heuristic method for solving
this latter non-convex optimization problem using the following ideas:

\begin{enumerate}

\item  We use a continuous relaxation of the non-convex, combinatorial problem
that arises by replacing the discontinuous step-function by a smoothed proxy
with parameter $\gamma$ allowing interpolation between the tractable convex
optimization problem (large $\gamma$) and the intractable non-convex,
combinatorial optimization ($\gamma = 0$).

\item We use the majorization-minimization algorithm to heuristically solve the
necessary non-convex optimization steps of this procedure by iteratively
linearizing the (concave) smoothed step function.

\end{enumerate}

\item Lastly, we extend all these methods by designing networks that are robust
against the failures of a small number of lines and/or generators. Essentially,
this is done by replacing the convex power-dissipation metric by the worst-case
power dissipation after removing some $k$ lines and/or generators.

\end{itemize}

The paper is structured as follows: (Section II-A) reviews the resistive network
model; (II-B) discusses how AC power flow is modeled by DC resistive network;
(III) presents the convex network optimization problem; (IV) presents the
non-convex extension to enforce sparsity; (V) presents robust network design;
(VI) indicates a number of potential extensions of our method and other
challenging open questions.

\section{Technical Preliminaries}

The optimization approach developed in Section III is based on the
resistive network model explained in Section II-A. We also
describe (Section II-B) that a modification of the effective resistive network is
adequate for the standard AC power flow model when considered
in the leading DC approximation.

\subsection{Resistive Network Model}

We give a brief introduction to electrical networks \cite{bollobas98,ds84}.
Let $G$ denote a graph with node set $N = \{1,\dots,n\}$ and $m$ (undirected) edges
$\{i,j\} \in G \subset 2^N$.  We assign edge weights $\theta_{ij} \equiv
\theta_{ji} \ge 0$ for all $\{i,j\} \in G$ ($\theta_{ij} = 0$ for
all non-edges $\{i,j\} \not\in G$). Regarded as a resistive network, the edges
$\{i,j\} \in G$ represent the \emph{lines} of the network  with $\theta_{ij}$
being the \emph{conductance} (inverse resistance) of a line.  We also use
$\ell \in G$ to index lines of the network. We define the \emph{conductance
matrix} $K(\theta) \in \mathbb{R}^{n \times n}$ of the network by  $K(\theta) =
\sum_{\{i,j\} \in G} \theta_{ij} (e_i - e_j) (e_i - e_j)^T$ where $e_i \in
\mathbb{R}^n$ are the standard basis vectors.  This is the edge-weighted
graph Laplacian of $G$ based on line conductances.  Thus,
\begin{equation}
K_{ij}(\theta) =
\left\{
\begin{array}{ll}
-\theta_{ij}, & i \neq j \\
\sum_{k \neq i} \theta_{ik}, & i = j
\end{array}
\right.,
\label{eq1}
\end{equation}
One may also write $K(\theta) = A \mathrm{Diag}(\theta) A^T$ where 
$\mathrm{Diag}(\theta) = \sum_\ell \theta_\ell e_\ell e_\ell^T \in \mathbb{R}^{m \times m}$ 
is a diagonal matrix and $A \in \mathbb{R}^{n \times m}$ is the incidence matrix of $G$ with 
columns $a_\ell = \pm (e_i - e_j)$ for each edge $\ell = \{i,j\}$.

Let $b \in \mathbb{R}^N$ represent the vector of \emph{injected
currents} --- nodes with $b_i > 0$ are sources, those with $b_i < 0$ are sinks
and $b_i = 0$ for transmission nodes.  In the resistive network, these represent
currents being injected into (or drawn from) each node by an external source.
Given $K$ and $b$, we obtain the (relative) electrical potential among the nodes
$u \in \mathbb{R}^n$ by solving the linear system of equations:
\begin{equation} \label{eq:kirchoff}
K u = b
\end{equation}
We observe the following properties of the conductance matrix (assuming connected $G$ and
non-zero $\theta$):

\begin{itemize}

\item $K$ is a symmetric positive semi-definite matrix: $u^T K u \ge 0$ for all
$u \in \mathbb{R}^n$.  As we will see later, this represents the fact that power
dissipation is non-negative.

\item $K$ has a single zero eigenvalue associated to the ``ones'' eigenvector:
$K \mathbf{1} = \mathbf{0}$.  This indicates that for $b = 0$ we must have 
uniform electric potential.

\item For any other eigenvector $K u = \lambda u$ (besides $u=\ones$) it holds that $\ones^T u = 0$
and $\lambda > 0$.

\end{itemize}
It is required that the total injected current is zero, $\ones^T b = \sum_{i \in
N} b_i = 0$, so that (\ref{eq:kirchoff}) can be satisfied.   Then, there is a
one-dimensional space of solutions of the form $\{u' + c \ones | c \in
\mathbb{R} \}$ for any $u'$ solving $K u' = b$, that is, the solution is
uniquely determined up to an overall additive shift of the electric potentials. 
There are several approaches one might use to ``regularize'' the problem of
computing $u$ such that the solution becomes unique.  Here, we require that
$\sum_i u_i = 0$,  obtained by solving the $n \times n$ system of equations $K'
u = b$ based on the invertible matrix $K' = K + \mathbf{1} \mathbf{1}^T$. One
may check that $K' \ones = n \ones$ and all other eigenvalues and eigenvectors
of $K'$ are the same as for $K$. The regularized solution to (\ref{eq:kirchoff})
is then given by $u = K^+ b$ where $K^+ \triangleq (K+\ones \ones^T)^{-1}$.

Current flow within the network is then determined by the electric potential $u$
and Ohm's law: the current flow from $i$ to $j$ is $b_{ij} = \theta_{ij} (u_i -
u_j)$.  Since $\theta_{ij} = \theta_{ji}$, it of course holds that $b_{ji} =
-b_{ij}$.  One may verify that $b_i + \sum_{k \neq i} b_{ki} = 0$ for all $i$
(current is conserved at each node). The total power loss over the network (due
to resistive heating of the lines) is given by:
\begin{equation}
\mathcal{L} =   \sum_{ij \in G} \theta_{ij} (u_i - u_j)^2 = u^T K u
\label{eq:Loss_Resist}
\end{equation}
Substitution of $u = K^+ b$ into this equation gives $\mathcal{L} = b^T K^+ K
K^+ b = b^T K^+ b$. If we fix the graph structure $G$ and the loads $b$, then
the power loss becomes a function of the conductances $\mathcal{L}(\theta) = b^T
(K(\theta)+\mathbf{1}\mathbf{1}^T)^{-1} b$. It is simple to generalize the power
loss objective to account for random fluctuations of the load $b$. For a random
current the expected power loss is:
\begin{eqnarray}
\mathcal{L}(\theta)
&=& \langle b^T K^+(\theta) b \rangle \nonumber \\
&=& \langle \tr(K^+(\theta) bb^T) \rangle \nonumber \\
&=& \tr(K^+(\theta) \langle bb^T \rangle) \nonumber \\
&=& \tr(K^+(\theta) B) \label{eq:KB}
\end{eqnarray}
where we have defined the matrix $B \triangleq \langle b b^T \rangle$, which is
a sufficient statistic of the random current model for the purpose of computing
the expected power loss. Importantly, $\mathcal{L}(\theta)$ is a convex
function, which is the basis for convex network optimization methods
\cite{gbs08,bvg01}.

\subsection{DC Approximation to AC Power Flow}

The existing power grid uses the AC voltages and currents generally described in
terms of complex amplitudes and lines with complex impedances, in contrast to
real currents and positive { conductances} of the resistive network setting. In
spite of this difference, the resistive network framework can be used to
approximate the AC system \cite{96WW}.

Indeed, (\ref{eq:Loss_Resist}) still holds in the case of AC flows if
$(u_i-u_j)^2$ is replaced by $|U_i-U_j|^2$, where $U_j$ is a complex potential
at the node $j$ of the network and $K$ now stands for the real part of the
network admittance matrix, also called the network (AC) conductance matrix
\cite{96WW}. In a healthy AC flow the voltage magnitude is stabilized to a
constant (unity in the rescaled  power units). In the so-called DC
approximation, where this stabilization is assumed ideal, $U_j=\exp({\it
i}\varphi_j)$ where real $\varphi_j$ is the phase of the potential and ${\it
i}^2=-1$. Susceptance of a transmission power line, defined as the imaginary
part of the line admittance, is normally an order of magnitude larger than the
respective real part (conductance of the line). Then the DC-approximation of the
AC Kirchhoff equations, with the conductance completely ignored, becomes
\begin{eqnarray}
\label{eq:p=betau}
p=\tilde{K} \varphi,
\end{eqnarray}
where $p$ is the vector of real power (with its components being
production/consumption at the graph nodes), and $\tilde{K}$ is the imaginary
part of the network admittance matrix, also called the network susceptance matrix,
and (\ref{eq:p=betau}) thus accounts only for the lossless transfer of real
power over the network, $\sum_i p_i=0$. We note that $\tilde{K}$ has all the
same essential properties of $K$ listed in Section II-A.

Substituting (\ref{eq:p=betau}) into the aforementioned expression for the power
losses over the network and keeping only the leading  DC-approximation terms
(first order in the conductance-to-susceptance ratio) one arrives at an
expression for losses
\begin{equation}
\mathcal{L} = \tfrac{1}{2} p^T \tilde{K}^+K\tilde{K}^+p,
\label{Loss_Resist}
\end{equation}
where $\tilde{K}^+ \triangleq (\tilde{K} + \ones \ones^T)^{-1}$. We
assume that the conductance-to-admittance ratio, $\mu$, is kept constant for all the
lines, i.e. $\tilde{K}=(1/\mu) K$. Then,  the only difference between the
DC-approximation model and the basic resistive network model will consist in
this additional re-scaling factor whose particular value is any case irrelevant
to the network optimization discussed in Section III. In particular,  this
translation from the resistive network model to the DC-approximation of the
AC-flow model means that (\ref{eq:KB}) turns into $\mathcal{L}(\theta)=
\tfrac{\mu^2}{2} \tr(K^+ B)$ where $B$ characterizes the random (real) power
flow through the network.

The main conclusion of this subsection is that with proper (and trivial)
rescaling the resistive network model is completely adequate to describe losses
in the leading order DC-approximation of the AC-flow model of the power grid.
Therefore,  with the understanding that we have neglected reactive power flows,
we may without loss of generality work with the resistive network model in the
remainder of the paper.

\section{Convex Network Optimization}

In this section we develop the main convex optimization method
we use to design electric power transmission networks.
This involves optimizing the line conductances for a given graph
to minimize the expected power loss subject to a linear constraint
(alternatively, adding a linear penalty) on the vector of line
conductances.  This is a generalization of the convex
optimization problem posed in \cite{gbs08}, which inspired our
approach of this paper.  The main contribution of this
section is in adapting their problem formulation to design electric
power transmission networks. In later sections, we also use
this convex optimization method as the core engine
within an iterative method for performing non-convex network
optimization with the aim of discovering good sparse network structures.

\subsection{The Network Optimization Problem}

First, we state the general form of the convex optimization problem that
we consider, and provide further details in the following subsections.
As discussed in Section II-A, we are given a graph $G$ of $n$ nodes
and $m$ edges.  The statistics of currents (power flows in the
DC-approximation of AC system) through the network are described by
an $n \times n$  matrix $B$.    Our aim is to assign the line conductances
$\theta$ to
balance the competing objectives of (1) maximizing network efficiency
(minimizing the expected power dissipation within the network) and (2)
minimizing the cost of building the network with conductances $\theta$.

We now specify a simple linear cost model on the line conductances. We model the
cost (say, in dollars) of building the network as $\alpha^T \theta = \sum_\ell
\alpha_\ell \theta_\ell$. The coefficients of this cost objective may be set as
$\alpha_\ell = c g^{-1} s_\ell^2$ where $c$ is the price of copper (per unit
volume), $g$ is the conductivity of copper and $s_\ell$ is the total length of
line $\ell$.  Then, $\alpha^T \theta$ represents the total cost of copper needed
to build the network with topology $G$, lines of length $s_\ell$ and
conductances $\theta_\ell$.  This follows as the conductance of a line of length
$s_\ell$ and cross-sectional area $a_\ell$ is $\theta_\ell = g a_\ell
s_\ell^{-1}$. Hence, the volume of a line is $s_\ell
a_\ell = s_\ell (g^{-1} s_\ell \theta_\ell) = g^{-1} s_\ell^2 \theta_\ell$ and
the cost of a line is $c g^{-1} s_\ell^2 \theta_\ell = \alpha_\ell \theta_\ell$.
Note that the problem of optimizing line conductances is essentially the same as
line sizing due to the linear correspondence between conductance and
cross-sectional area.

Given $G$, $B$ and $\alpha$ one may then select the line conductances $\theta$
to make the network as efficient as possible (minimizing the expected power loss
due to resistive heating of the lines) subject to a linear constraint that the
total cost of building the network must be no greater than a specified budget
$C$:
\begin{displaymath}
\begin{array}{ll}
\mbox{minimize} & \mathcal{L}(\theta) \\
\mbox{subject to} & \theta \ge 0 \\
& \alpha^T \theta \le C
\end{array}
\end{displaymath}
This is essentially the same as the convex optimization problem posed in
\cite{gbs08}. The \emph{total resistance} metric that they considered is
recovered by setting $B$ equal to the identity matrix.   This was interpreted as
the expected power loss under a Gaussian random current model $b \sim
\mathcal{N}(0,I)$ (modulo a projection to enforce the constraint $\ones^T b =
0$). Equivalently, one may replace the budget constraint by a linear penalty on
network cost, solving the convex optimization problem:
\begin{displaymath}
\min_{\theta \ge 0} \left\{ \mathcal{L}(\theta) + \lambda \alpha^T \theta
\right\}
\end{displaymath}
The parameter $\lambda > 0$ is a Lagrange multiplier enforcing the budget
constraint (the two problems are equivalent for corresponding values of $C$ and
$\lambda$). Alternatively, we may set
$\lambda^{-1} = p T$ where $p$ is the cost of power generation and $T$ is the
expected operational lifetime of the network. Then, the solution of the
penalized optimization problem yields the most cost-effective network design,
minimizing the sum of the cost to build the network and the cost to operate the
network over its operational lifetime.   In the remainder of the paper, we focus
of this latter form of the problem setting $\lambda = 1$ (redefining $\alpha
\rightarrow \lambda \alpha$).

Our main contribution in the remainder of the section is to further tailor this
problem to the setting of electric power transmission by appropriate definition
of $B$.

\subsection{Single-Generator Formulation}

First, we address the simplest case of optimizing a network
with multiple independent random loads supplied by a single generator
at a specified location. For non-generator nodes we specify
the mean load $\bar{b}_i = \langle b_i \rangle < 0$ and
the standard deviation $\sigma_i = \langle (b_i-\bar{b}_i)^2
\rangle^{\frac{1}{2}}$.
At transmission nodes we set $\bar{b}_i = 0$ and $\sigma_i = 0$.
At the generator node we must have $b_0 = -\sum_{i \neq 0} b_i$ to satisfy
the constraint $\sum_i b_i = 0$. Then, the overall random load matrix $B =
\langle b b^T \rangle$
is given by:
\begin{displaymath}
B =
\left(
\begin{array}{cc}
(\sum_{i \neq 0} \bar{b}_i)^2 + \sum_{i \neq 0} \sigma_i^2 & - \ones^T (\bar{b}
\bar{b}^T + \Sigma) \\
-(\bar{b} \bar{b}^T + \Sigma) \ones & \bar{b} \bar{b}^T + \Sigma \\
\end{array}
\right)
\end{displaymath}
where $\Sigma = \mathrm{Diag}(\sigma)$ is the diagonal covariance
matrix of non-generator loads.  One could also use a general
covariance matrix $\Sigma$ if cross-correlations among the consumers
is known (e.g., induced by hidden variables such as the time, season or
environmental factors such as temperature).

\subsection{Multiple-Generator Formulation}

Next, we consider the case that there are two or more generators
within the network.  One could consider explicitly modeling the
full matrix $B$, including both power consumption and generation.
However, for controlled power generators this is not realistic
because the response of generators to meet demand will surely depend
on the network itself (being designed) and moreover is adaptive
to fluctuations in the spatial distribution of demand.   To provide
a more realistic model of power generation, we will assume that
power generation is always chosen optimally in response to demand
and network configuration.  That is to say, for any given demand
$b_c$ the generation $b_g$ is chosen subject to minimize
$b^T K^+(\theta) b$ where $b = (b_c, b_g)$, subject to the constraint
$\ones^T b_g = - \ones^T b_c$.  Then, averaging the optimized
power loss over the distribution of $b_c$ leads to a new
convex objective $\hat{L}(\theta)$ that we may use in
the convex network optimization problem.

Although this may at first appear to be more complicated,  it
turns out there is a simple trick that allows us to transform
it back to the problem we have already considered.
Let $G'$ be an augmented representation of the network
in which we include one auxiliary node $0$, considered as a virtual
generator, and where we add auxiliary lines connecting this virtual
generator to each of the real generator nodes of $G$.  Now,
we may apply the formulation of Section III-B to this augmented
model, where the virtual generator is treated as the only generator
and the actual generator nodes of $G$ are now treated simply as
transmission nodes.  By setting the conductance of virtual lines
to infinity, we can allow current (power) to flow freely without
dissipation from the virtual generator to the real generators.
Thus, solving for power flows in this augmented network model
uses the optimal flow (minimizing power dissipation) and is
equivalent to optimizing the power generation in the original
model.  We omit technical proofs, which essentially involve showing that
the current flow described by Kirchoff's laws is efficient.

Finally, rather then actually setting the virtual lines to have
infinite conductance, we can make the cost of conductance
on these lines negligible in comparison to real lines, so that
the virtual lines are assigned very large conductances (relative
to real lines) in the solution of the convex network optimization
problem.

\subsection{Convex Optimization Algorithm}

In this section we briefly describe the method we use
to solve the convex network optimization problem.  The
main technical result ones needs are formulas for
the gradient vector and Hessian matrix of the expected
power loss $\mathcal{L}(\theta)$.
Generalizing those derivations of \cite{gbs08}, one
obtains:
\begin{displaymath}
\nabla \mathcal{L}(\theta) = -\frac{1}{2} \mathrm{diag}(A^T K^+(\theta) B
K^+(\theta) A)
\end{displaymath}
\begin{displaymath}
\nabla^2 \mathcal{L}(\theta) = (A^T K^+(\theta) A) \circ (A^T K^+(\theta) B
K^+(\theta) A)
\end{displaymath}
Similar to \cite{bvg01,gbs08}, we enforce the non-negativity constraint $\theta
\ge 0$
using the log-barrier function \cite{bv04}:
\begin{displaymath}
\min_{\theta > 0} \left\{ \mathcal{L}(\theta) + \alpha^T \theta - \zeta
\sum_{\ell \in G} \log \theta_\ell \right\}
\end{displaymath}
The solution of this modified problem will always be strictly positive.
One may obtain a close approximation to the optimal solution of the original
problem
for sufficiently small values of $\zeta$ \cite{bvg01}.  Efficient algorithms
start
by (approximately) solving this problem for a large value of $\zeta$ and then
iteratively updating the solution for a decreasing sequence of $\zeta$ values.
It is straight-forward, using the formula above, to implement Newton's method
with back-tracking line search to minimize this convex objective function
\cite{bv04}.

We remark one technical difficulty we have encountered. Using the formula
$\mathcal{L}(\theta) = \tr(K^+(\theta) B)$, it may not always be possible to
make $\zeta$ arbitrarily small.  This is due to numerical difficulties with
computing $\mathcal{L}(\theta)$ when the graph $G$ is becoming effectively
disconnected due to many $\theta$'s going to zero.  The matrix $K(\theta)+\ones
\ones^T$ is becoming singular in such cases, such that the formula for
$\mathcal{L}(\theta)$ should really be reformulated with respect to a subgraph
of $G$ with non-zero conductances. However, these technical difficulties may be
avoided by not letting $\zeta$ become so small that $K(\theta)+\ones \ones^T$
becomes numerically singular.  In future work, it may be desirable to develop a
robust way of computing $\mathcal{L}(\theta)$ in such cases so that $\zeta$ can
be made arbitrarily small.

\subsection{Demonstrations}

We now describe our fist set of demonstrations, based on four examples that we
revisit in later sections.   All examples having essentially the same graph
topology but with different configurations of demand and generation nodes.  The
graph $G$ is comprised of a $w \times w$ grid of nodes (with $w = 9$ or $10$)
and has lines  between nearest and second-nearest neighbors of the grid
(resulting in vertical/horizontal edges between nearest neighbors and diagonal
edges between second-nearest neighbors).  Transmission, consumption and
generation nodes are respectively marked as black dots, blue dots and red dots.
  We set $\alpha = 1$ on horizontal/vertical edges and $\alpha = 2$ on diagonal
edges. We have $\bar{b} = -1$ and $\sigma = \tfrac{1}{3}$ at consumer nodes; and
$\bar{b} = \sigma = 0$ at transmission nodes.  Fig.~1 shows the result of
solving the convex network optimization in four examples. In the multiple
generator case (seen at lower-right), we do not show the virtual generator
or the lines to this generator.  We observe that:
\begin{itemize}

\item The solution is somewhat sparse in these examples (it does not use all
edges of $G$) but is not as sparse as possible (it is not a minimal tree/forest
needed to connect consumers to generators).

\item It is not necessary that all transmission nodes are involved in the
solution, as shown by the example seen at the lower-left of the figure.

\item With multiple generators (lower-right), the power transmission network may
become disconnected, with each generator serving a particular subset of nearby
consumer nodes.

\end{itemize}

\begin{figure}

\noindent \hspace{-.6cm} \includegraphics[scale=.113]{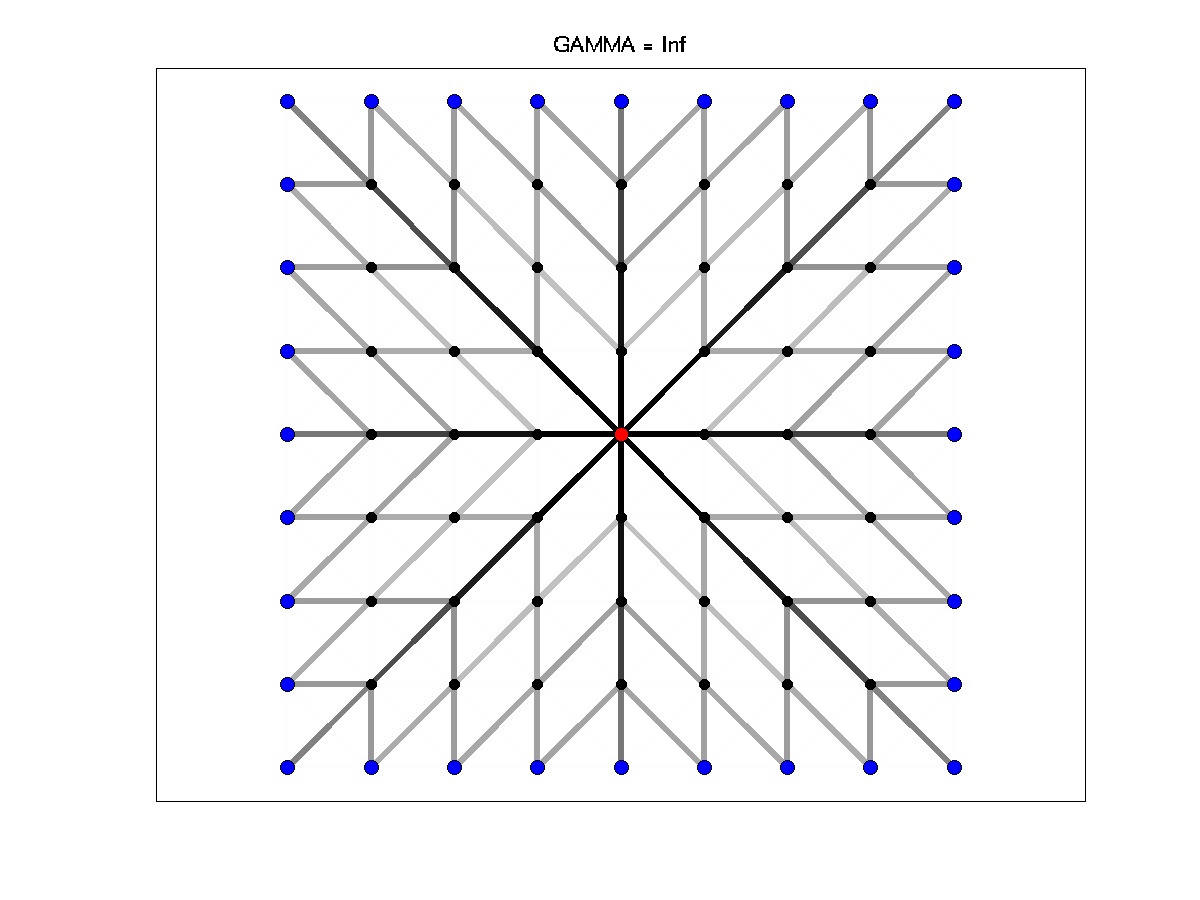}
\hspace{-.6cm} \includegraphics[scale=.113]{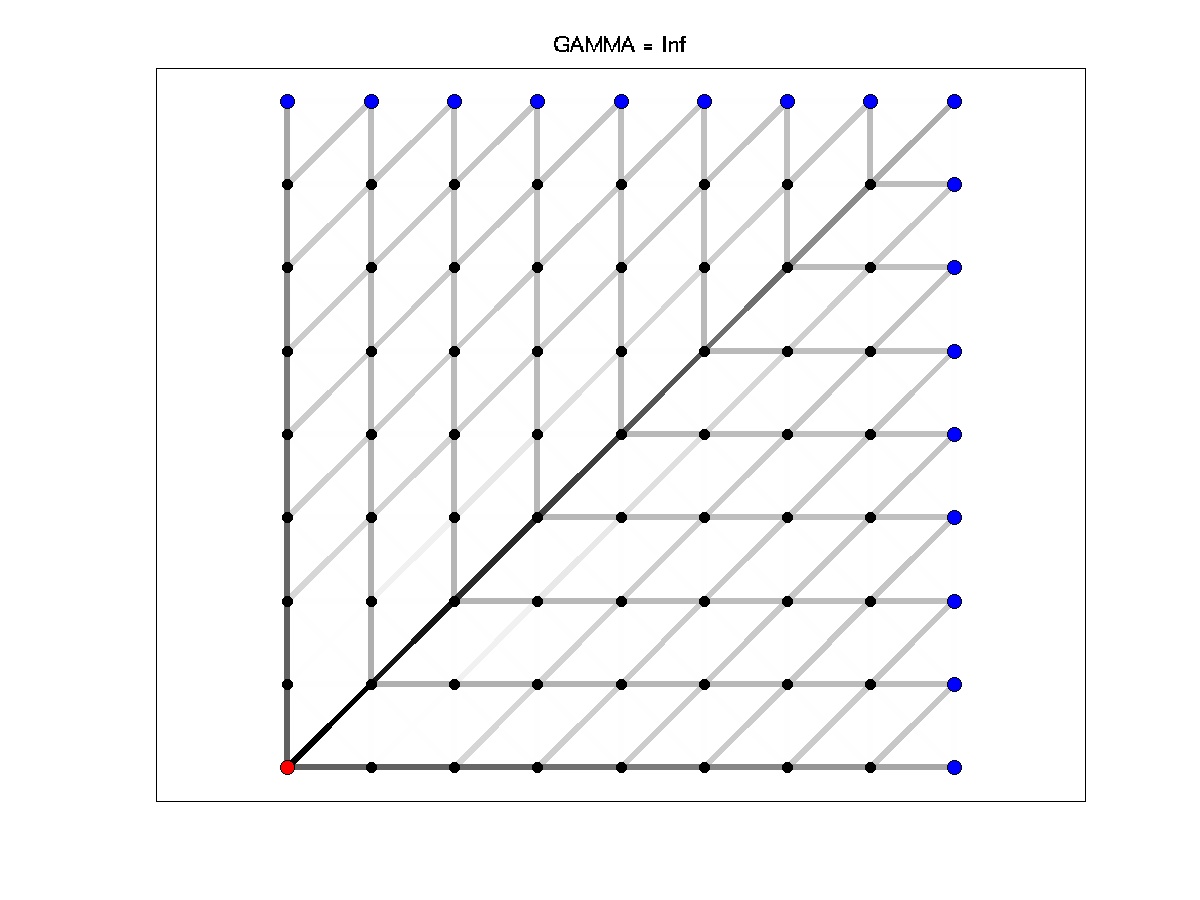}

\vspace{-.2cm}
\noindent \hspace{-.6cm} \includegraphics[scale=.113]{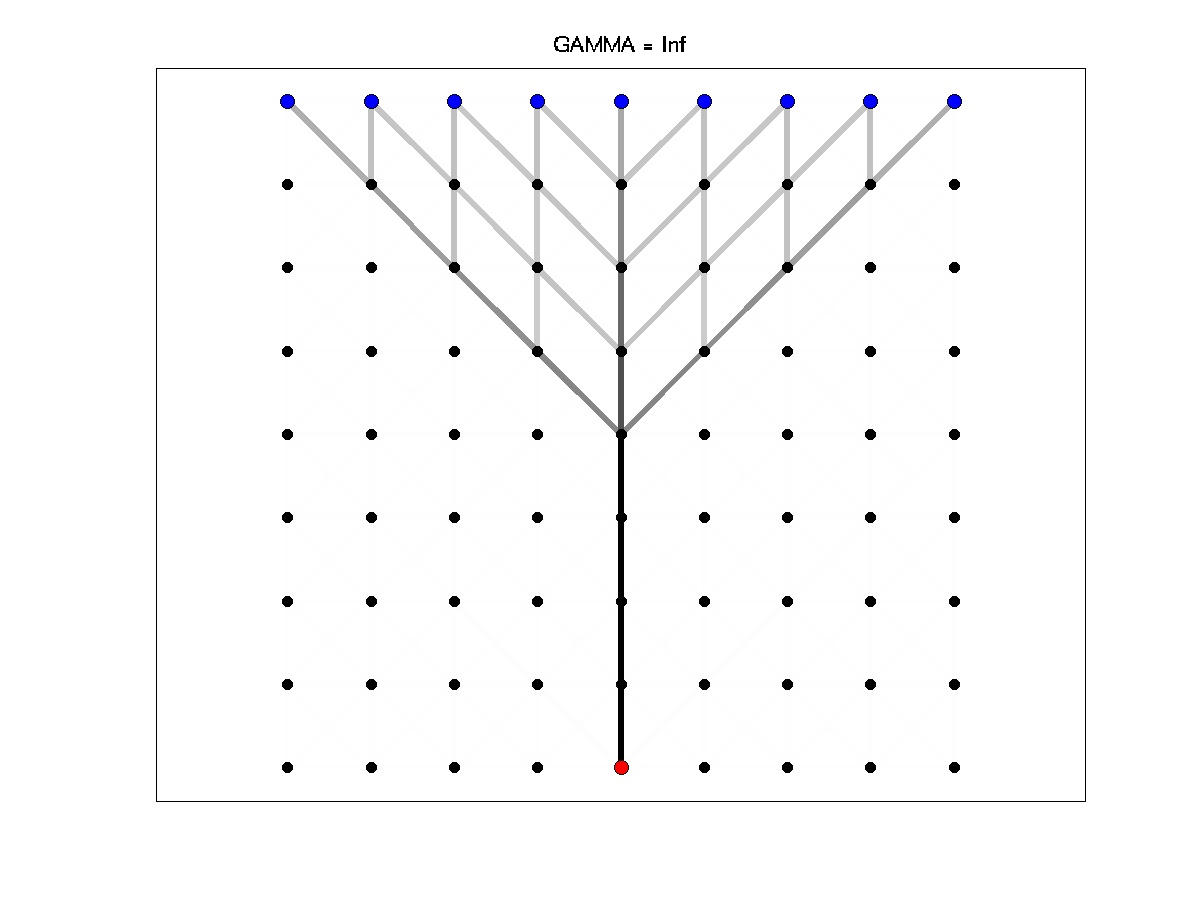}
\hspace{-.6cm} \includegraphics[scale=.113]{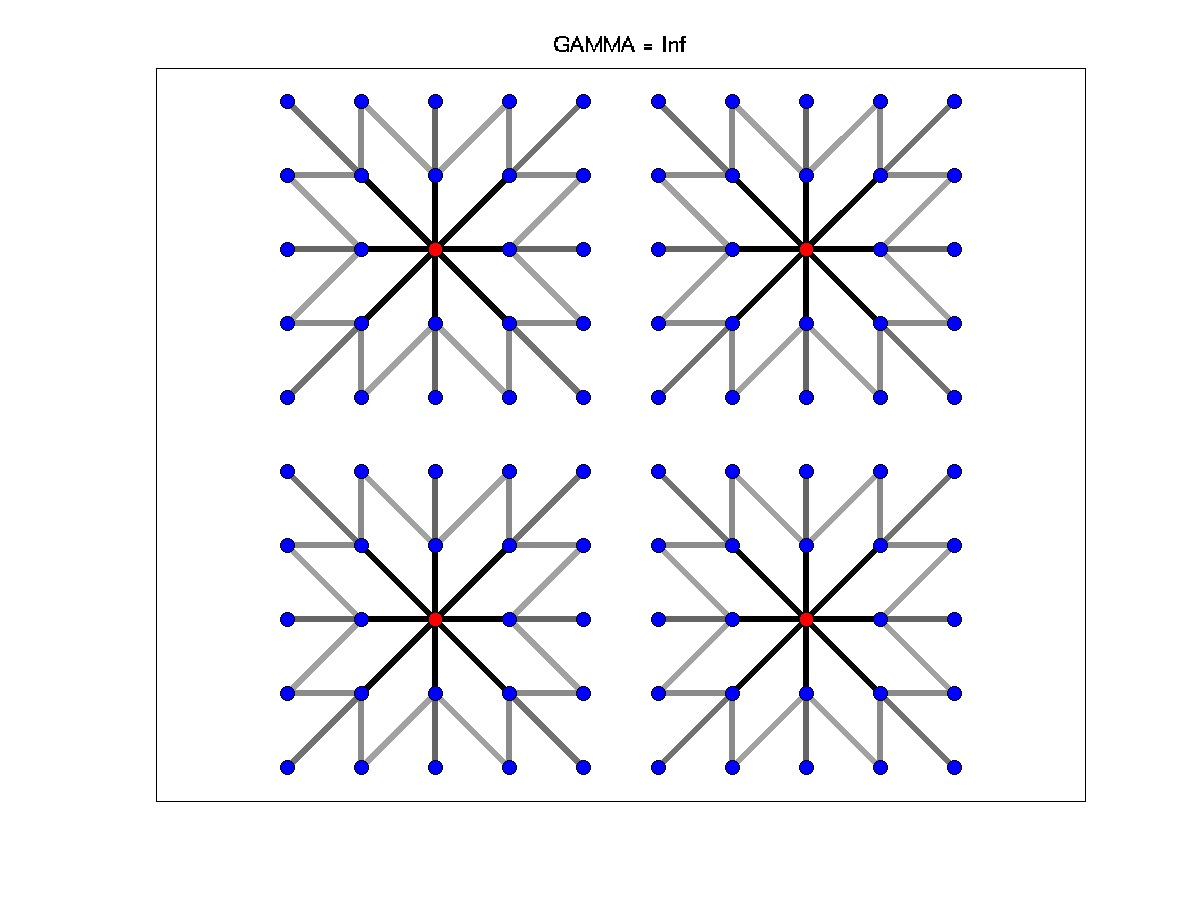}

\vspace{-.3cm}
\caption{Illustration of globally optimal network designs in the convex
network optimization method.  The strength of a line (conductance) is
indicated by the darkness of the drawn edge, such that zero-conductance
lines are not seen.}
\vspace{-.3cm}

\end{figure}

\section{Selecting Network Structure}

In this section we present a non-convex generalization of the approach taken in
the preceding section.  As we have seen, the convex optimization method does not
always produce sparse solutions, i.e., it will typically use most of the edges
of the graph $G$.  In practical applications, we expect that such solutions are
undesirable, as we would like to use the simplest network (with as few edges as
possible) that is sufficient to meet power transmission requirements.  Towards
this end, we reformulate the network cost part of our optimization objective so
as to favor solutions with fewer edges.  However, this then gives a non-convex
optimization problem, which is generally intractable to solve exactly.  Hence,
we develop a heuristic approach to find good solutions of this non-convex
optimization problem.  Using the majorization-minimization algorithm, we are
able to approximately solve the non-convex problem by instead solving a sequence
of convex optimization problems.  Moreover, each convex optimization problem
will be of the form solved in Section III, with the vector $\alpha$ being
iteratively modified.  Thus, the methods of Section III provide an optimization
engine for the non-convex optimization method developed in this section.

\subsection{Sparsity-Favoring Network Cost}

In practice, we may also require that the network should by sparse.
We formulate this by adding a cost on all lines with non-zero
conductance so as to encourage solutions with as few lines as possible:
\begin{displaymath}
\min_{\theta \ge 0} \{ \mathcal{L}(\theta) + \alpha^T \theta + \beta^T
\phi(\theta) \}
\end{displaymath}
where $\phi(t)$ is the unit-step function, $\phi(0) = 0$ and $\phi(t)=1$ for all
$t>0$, which is applied element-wise to $\theta$ such that $\beta^T \phi(\theta)
= \sum_\ell \beta_\ell \phi(\theta_\ell)$.  Note the $\phi(t)$ is non-convex (in
fact, it is concave on $t \ge 0$) and discontinuous at $t=0$,

Whereas the linear cost $\alpha^T \theta$ essentially represents the cost of
copper needed to build the network, the non-convex cost $\beta^T \phi(\theta)$
represent other costs that do not scale with conductance, e.g., the cost of
purchasing/leasing right-of-way along corridors of land along the lines and
other expenses (labor, poles, towers, environmental impact) that do not tend to
zero for low-conductance lines.  More realistically, we expect the actual cost
of a line to be a concave function of conductance approaching an affine function
for large enough conductances.  The simple model above roughly captures this
behavior.

However, minimizing this discontinuous objective function is now a difficult
combinatorial optimization problem.  In order to find the optimal
solution, one would need to enumerate all possible subgraphs of $G$ and
perform a convex network optimization within each subgraph.  Clearly
this is not a scalable approach.  We instead propose a heuristic
solution method in the following subsections.

\subsection{Annealed Smoothing Method}

To avoid having to perform a combinatorial optimization, we begin
by smoothing the objective function to a continuous (albeit non-convex)
objective.  This is accomplished by replacing the step function $\phi(t)$
by a continuous penalty function with smoothing parameter $\gamma > 0$:
\begin{displaymath}
\phi_\gamma(t) = \frac{t}{t+\gamma}
\end{displaymath}

We observe that $0 \le \phi_\gamma(t) \le 1$ (for $t \ge 0$), $\phi_\gamma(0)=0$
and $\phi_\gamma(t) \rightarrow 1$ as $t \rightarrow \infty$.  Hence,
$\phi_\gamma$ may be regarded as a smooth approximation to the step function
$\phi$. Moreover, $\phi_\gamma(t) \rightarrow \phi(t)$ as $\gamma \rightarrow
0$, such that the optimal solution of the smoothed problem should approach that
of the combinatorial problem for sufficiently small $\gamma$.  However, because
this gives a non-convex optimization problem, there may be many local minima and
it can still be intractable to determine the global minimum.

Next, we observe that for very large values of $\gamma$ we have $\phi_\gamma(t)
\approx 0$ (over a large range of $t$) so that the smoothed problem becomes
equivalent to the convex optimization problem.  This suggests an ``deterministic
annealing'' strategy in which we start with the solution of the tractable convex
problem (corresponding to large $\gamma$) and then iteratively update the
solution while gradually decreasing $\gamma$ to a small value (approximating the
difficult combinatorial problem). At each stage of this annealing procedure, we
must solve a non-convex minimization starting from an initial guess
corresponding to a local minima of the preceding optimization.

\subsection{Majorization-Minimization Algorithm}

The approach described above requires solving (heuristically) a sequence of
non-convex optimization problems.  In this section, we present the algorithm
used for this purpose.  It is an instance of the general
majorization-minimization algorithm \cite{hl04}, which we now review.  Consider
an objective function of the form
\begin{displaymath}
f(x) = f_\cup(x) + f_\cap(x)
\end{displaymath}
where $f_\cup$ is a convex function and $f_\cap$ is a concave function.  In
order to seek a local minimum of $f$, one may iteratively linearize the concave
part $f_\cap$ and minimize the resulting convex upper-bound to $f$.  That is,
given the previous guess $x^{(k)}$ of the solution, we may then approximate
$f_\cap(x)$
by its linear upper-bound:
\begin{displaymath}
f_\cap(x) \le f_\cap(x^{(k)}) + \nabla f_\cap(x^{(k)})^T (x - x^{(k)}).
\end{displaymath}
This gives a convex upper-bound approximation of the objective:
\begin{displaymath}
f(x) \le f_\cup(x) + \nabla f_\cap(x^{(k)})^T x + \mathrm{const}
\end{displaymath}
where $\mathrm{const}$ is independent of $x$.  We then minimize this convex
function to obtain the next guess:
\begin{displaymath}
x^{(k+1)} = \arg\min_x \left\{ f_\cup(x) + \nabla f_\cap(x^{(k)})^T x \right\}.
\end{displaymath}
This procedure is guaranteed not to increase the objective $f(x)$ and typically
converges to a local minimum of $f(x)$. However, saddle-points and local maxima
are also unstable fixed-points of the algorithm.  Adding a small
random perturbation to $x$ will cause the method to leave such non-minimal
fixed-points.

It is straight-forward to apply this method in our present setting with $x =
\theta$,
$f_\cup(\theta) = \mathcal{L}(\theta) + \alpha^T \theta$ and $f_\cap(\theta) =
\beta^T \phi_\gamma(\theta)$
(note that $\phi_\gamma(t)$ is concave).  Doing so, we obtain the following
iterative algorithm:
\begin{eqnarray}
\alpha_\ell^{(k)} &=& \alpha_\ell + \tfrac{\gamma}{\left(\gamma +
\theta^{(k-1)}_\ell\right)^2} \beta_\ell \nonumber \\
\theta^{(k)} &=& \arg\min_{\theta \ge 0} \left\{ \mathcal{L}(\theta) +
(\alpha^{(k)})^T
\theta \right\}
\end{eqnarray}
Observe that the optimization problem is of the same form that we considered in
the Section III (with a modified value of $\alpha$) and can hence be
solved using the methods of that section.

To accelerate convergence of the convex optimization algorithm, the parameter
$\zeta$ of the barrier method may be kept fixed to a small value after the
initial optimization. This is usually more efficient because small changes of
the smoothing parameter $\gamma$ or the coefficient vector $\alpha$ typically
do not produce a large change in the optimal $\theta$.   However, occasionally
it can happen that reducing $\gamma$ can cause a local minima to disappear, such
that the algorithm must migrate to another local minima.  When this happens, as
indicated by Newton's method not converging after a reasonable number of
iterations, it is better to ``restart'' the barrier method with a large value of
$\zeta$.

\subsection{Demonstrations}

We return to those four examples introduced in Section III-E.  Setting $\beta =
1$ on horizontal/vertical edges and $\beta = \sqrt{2}$ on diagonal edges, we
obtain the solutions seen in Fig.~2 using the smoothed-annealing and
majorization-minimization methods described in this section.  In these examples,
the graph is ``thinned'' to a minimal tree or forest sufficient to distribute
power to consumer nodes from the generator(s).  However, the level of thinning
actually depends on the relative sizes of $\alpha$ and $\beta$, e.g., for
smaller values of $\beta$ one can obtain intermediate solutions between those
seen in Figs.~1 and 2.

\begin{figure}

\noindent \hspace{-.6cm} \includegraphics[scale=.113]{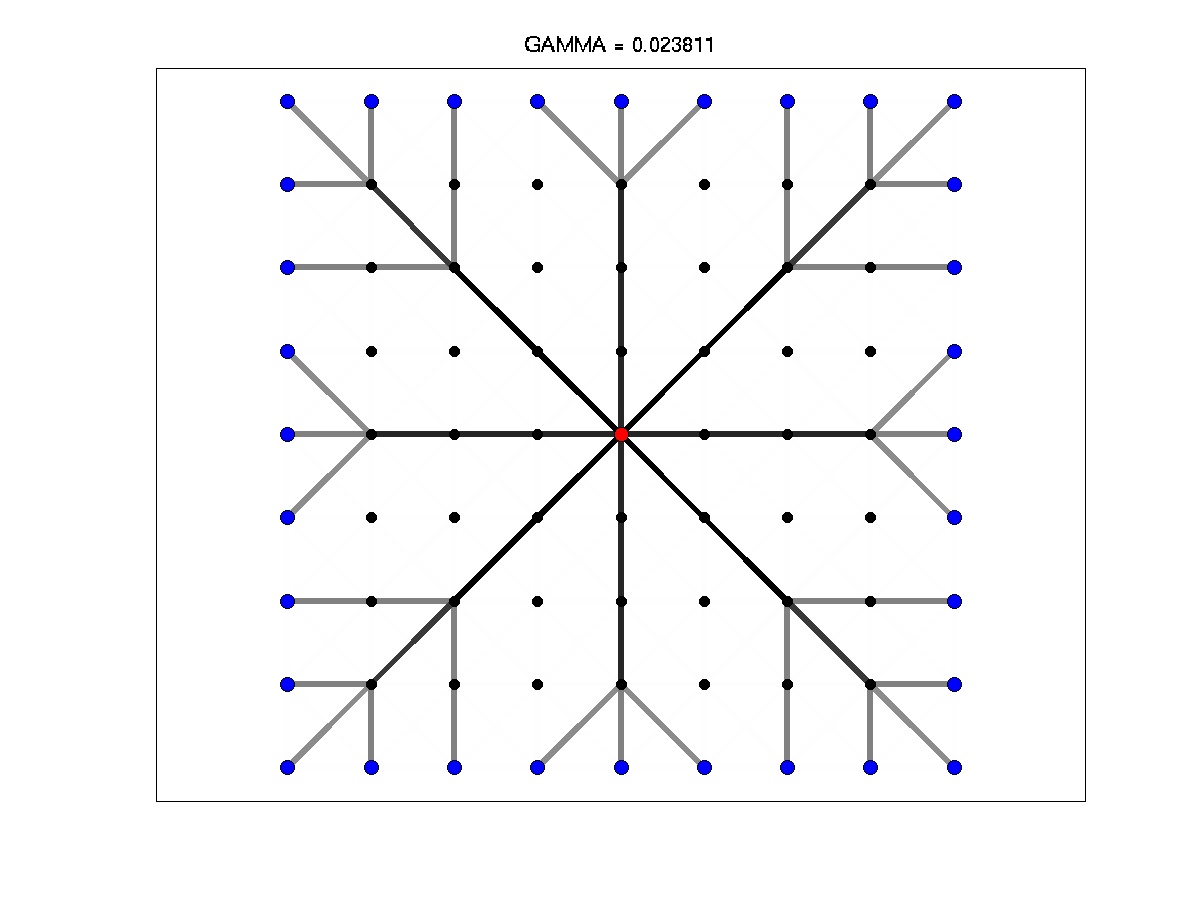}
\hspace{-.6cm} \includegraphics[scale=.113]{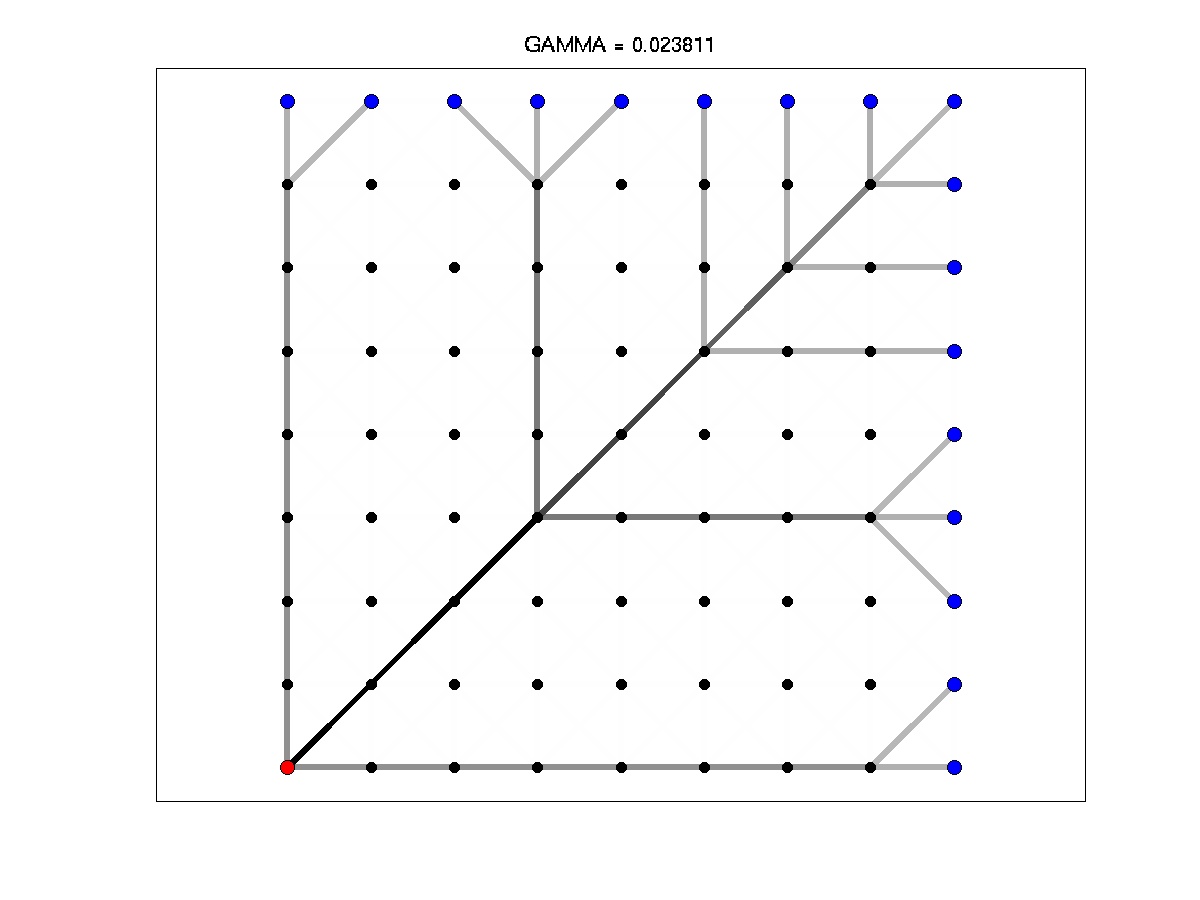}

\vspace{-.2cm}
\noindent \hspace{-.6cm} \includegraphics[scale=.113]{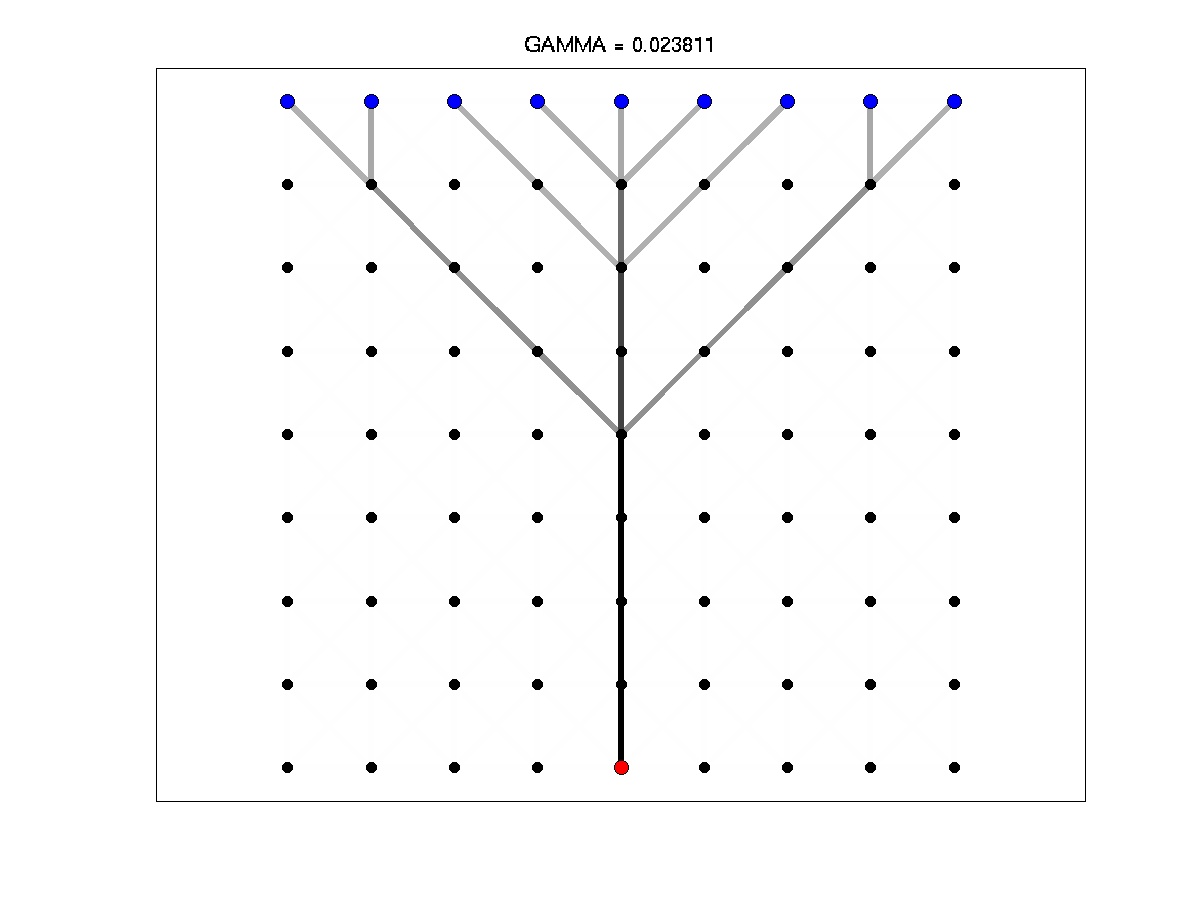}
\hspace{-.6cm} \includegraphics[scale=.113]{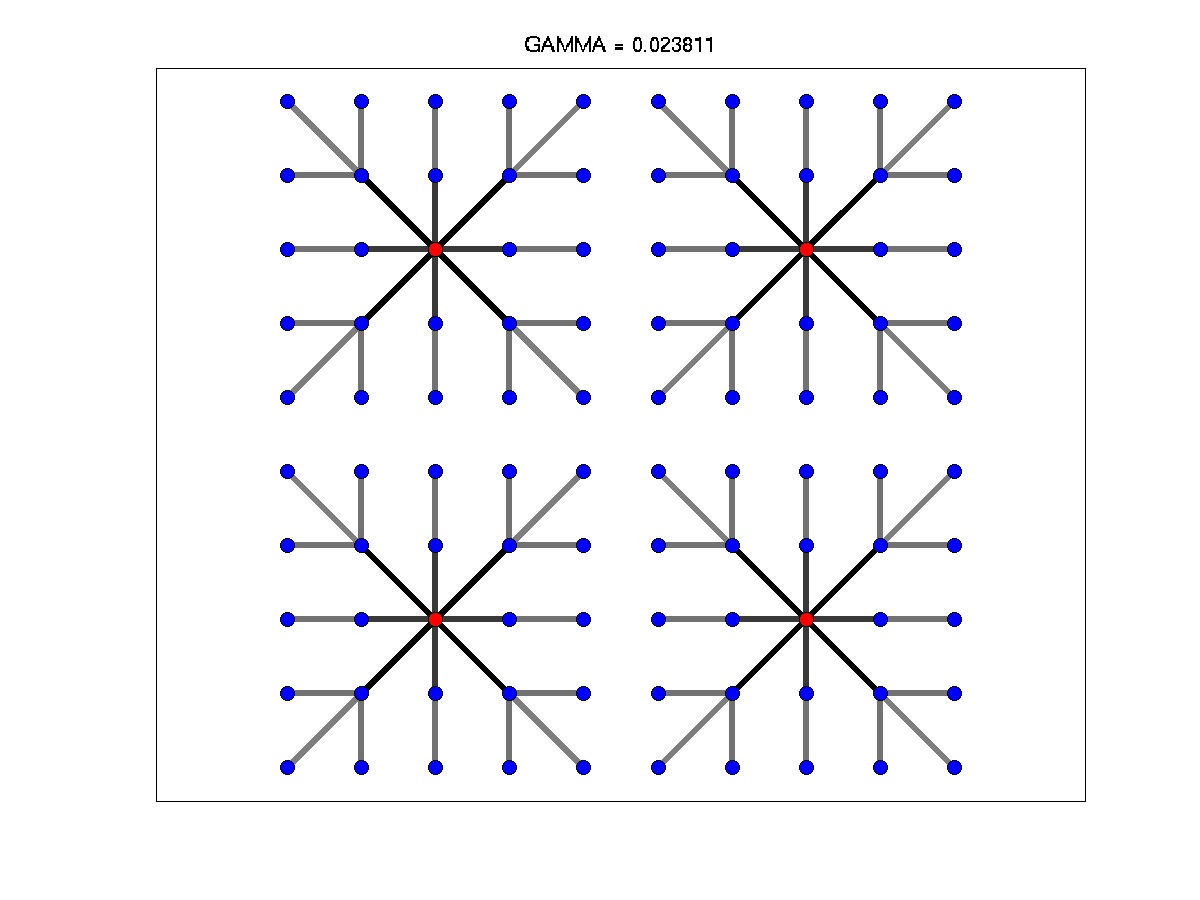}

\vspace{-.3cm}
\caption{Illustration of sparse network designs.  Compare to Fig.~1.}
\vspace{-.3cm}

\end{figure}

\section{Robust Network Design}

In this section we propose a simple modification of the methods of the previous
sections to obtain network designs that are robust to the failure of lines
and/or generator of the network.

\subsection{Imposing Robustness to Line/Generator Failures}

First, we observe that the expected power loss $\mathcal{L}(\theta)$
acts as a kind of barrier function which enforces the constraint that the
network must remaining connected such that every consumer node is connected to
at least one generator node.  For instance, if we gradually reduce $\theta$ to
zero on any subset of edges that separates a consumer node from all generators,
the power loss $\mathcal{L}(\theta)$ will tend to infinity.

Now, suppose that lines of the power network are subject to failures, meaning
that the conductance of a line is set to zero.  We would like for the network to
still be able to supply the consumers and do so without large power dissipations
(e.g., caused by having to route large currents through low-conductance lines
after a failure).  Let $z \in \{0,1\}^m$ be an indicator vector of line failures
such that $z_\ell = 1$ for failed lines and $z_\ell = 0$ for operational lines.
The power dissipation after removing failed lines is $\mathcal{L}(\theta; z)
\triangleq \mathcal{L}((\ones-z) \circ \theta)$. Suppose that we require that
the network is robust to up to $k$ line failures. Then, the worst case power
dissipation is:
\begin{displaymath}
\mathcal{L}^{\setminus k}(\theta) = \max_{\ones^T z = k} \mathcal{L}(\theta; z)
\end{displaymath}
Importantly, we note that this robust power dissipation is a convex function of
$\theta$ because the point-wise maximum of a collection of convex functions is
also a convex function.  We may then design our power network so as to minimize
the robust objective:
\begin{displaymath}
\min_{\theta \ge 0} \left\{ \mathcal{L}^{\setminus k}(\theta) + \alpha^T \theta
+ \beta^T \phi(\theta) \right\}
\end{displaymath}
Again using the annealed smoothing and majorization-minimization methods, this
results in a sequence of convex optimization problem of the form
$\min_\theta \left\{ \mathcal{L}^{\setminus k}(\theta) + \tilde{\alpha}^T
\theta\right\}$. We note that the robust power loss serves as a barrier function
to ensure that the graph remains $(k+1)$-connected, such that every consumer node
must be connected to the (virtual) generator by at least $k+1$ distinct paths.

Lastly, we remark that we may apply the same method in the multiple-generator
setting by requiring robustness with respect to failures of lines to the
virtual generator in the augmented network representation, which is equivalent
to allowing for failures of generators.

\subsection{Gibbsian ``Soft-Max'' Optimization}

While this approach results in a convex optimization problem, its solution using
standard steepest descent methods such as Newton's method is complicated by the
fact that $\mathcal{L}^{\setminus k}$ is non-smooth (it is not everywhere
differentiable). One could handle this using non-smooth optimization methods
such as subgradient descent.  However, these methods tend to converge slowly.
We avoid this complication by using the Gibbsian ``soft-max'' function with
smoothing parameter $\tau > 0$:
\begin{displaymath}
\mathcal{L}^{\setminus k}_\tau(\theta) = \tau \log \sum_{\ones^T z = k}
\exp\left[ \tau^{-1} \mathcal{L}(z ; \theta) \right]
\end{displaymath}
This is a smooth, convex function of $\theta$ and gives an upper-bound to
$\mathcal{L}^{\setminus k}(\theta)$. Moreover, $\mathcal{L}^{\setminus k}_\tau
\rightarrow \mathcal{L}^{\setminus k}$ (uniformly) as $\tau \rightarrow 0$.
Thus, smooth convex optimization of $\mathcal{L}^{\setminus k}(\theta)$ provides
a good approximation to the non-smooth convex optimization of
$\mathcal{L}^{\setminus k}$ for sufficiently small values of $\tau$.

\subsection{Demonstrations}

Again, we return to those four examples discussed previously in Sections III-E
and IV-D. Using the smoothed robust power loss with $k=1$ and $\tau=.01$, we
obtain the robust network solutions seen in Fig.~3 using the smoothed annealing
and majorization-minimization methods.   Observe that these graphs are
two-connected, such that every consumer node is connected to the (virtual) generator
nodes by at least two distinct paths.  Also, the multiple generator solution is now a
connected graph so as to be robust to failure of a generator.

\begin{figure}

\noindent \hspace{-.6cm} \includegraphics[scale=.113]{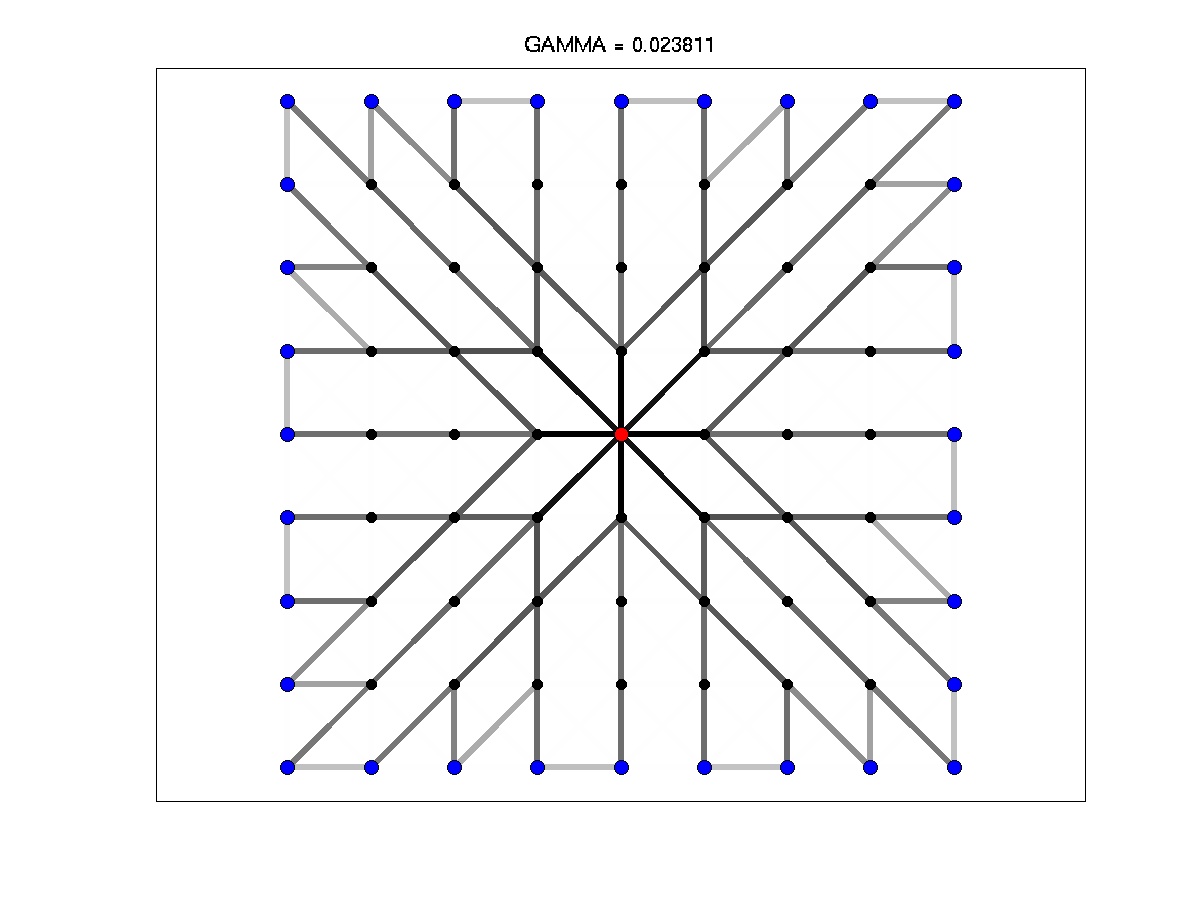}
\hspace{-.6cm} \includegraphics[scale=.113]{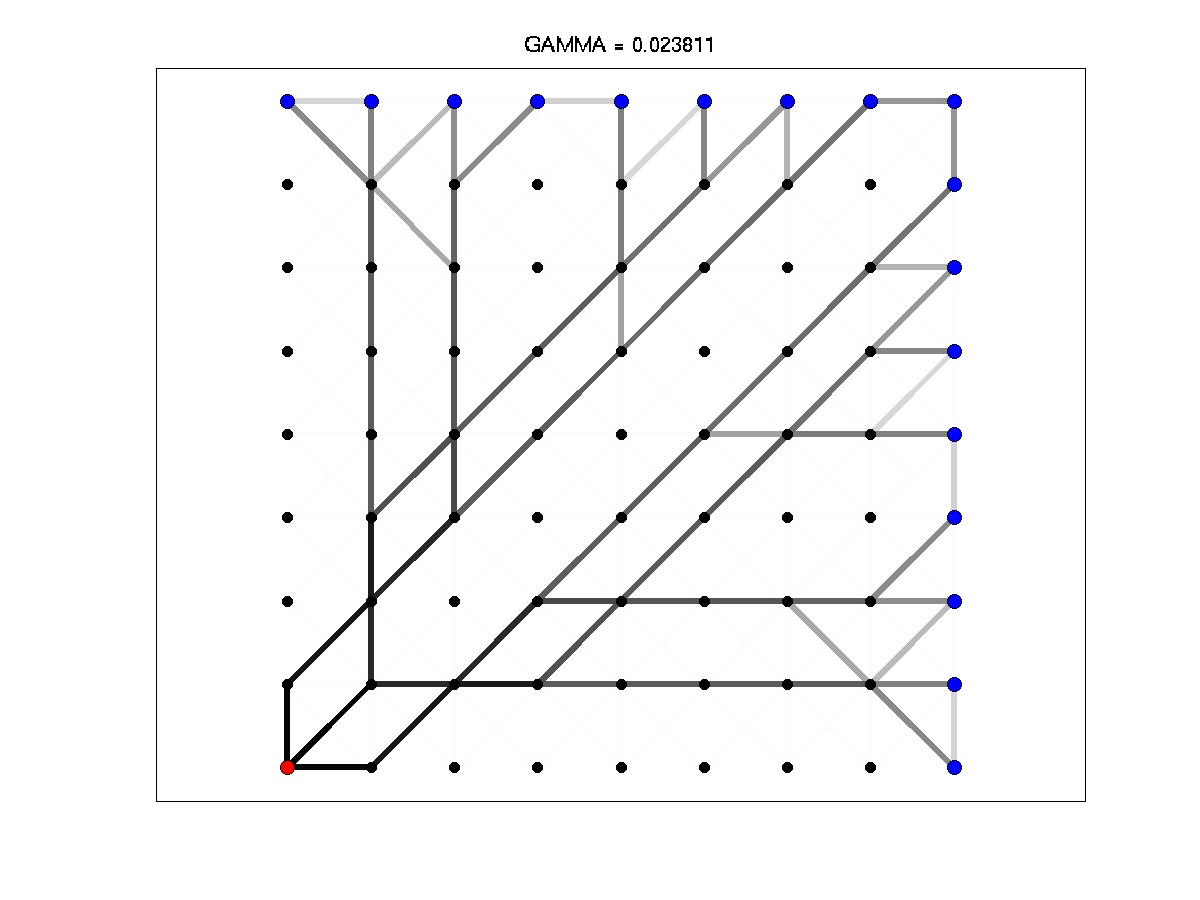}

\vspace{-.2cm}
\noindent \hspace{-.6cm} \includegraphics[scale=.113]{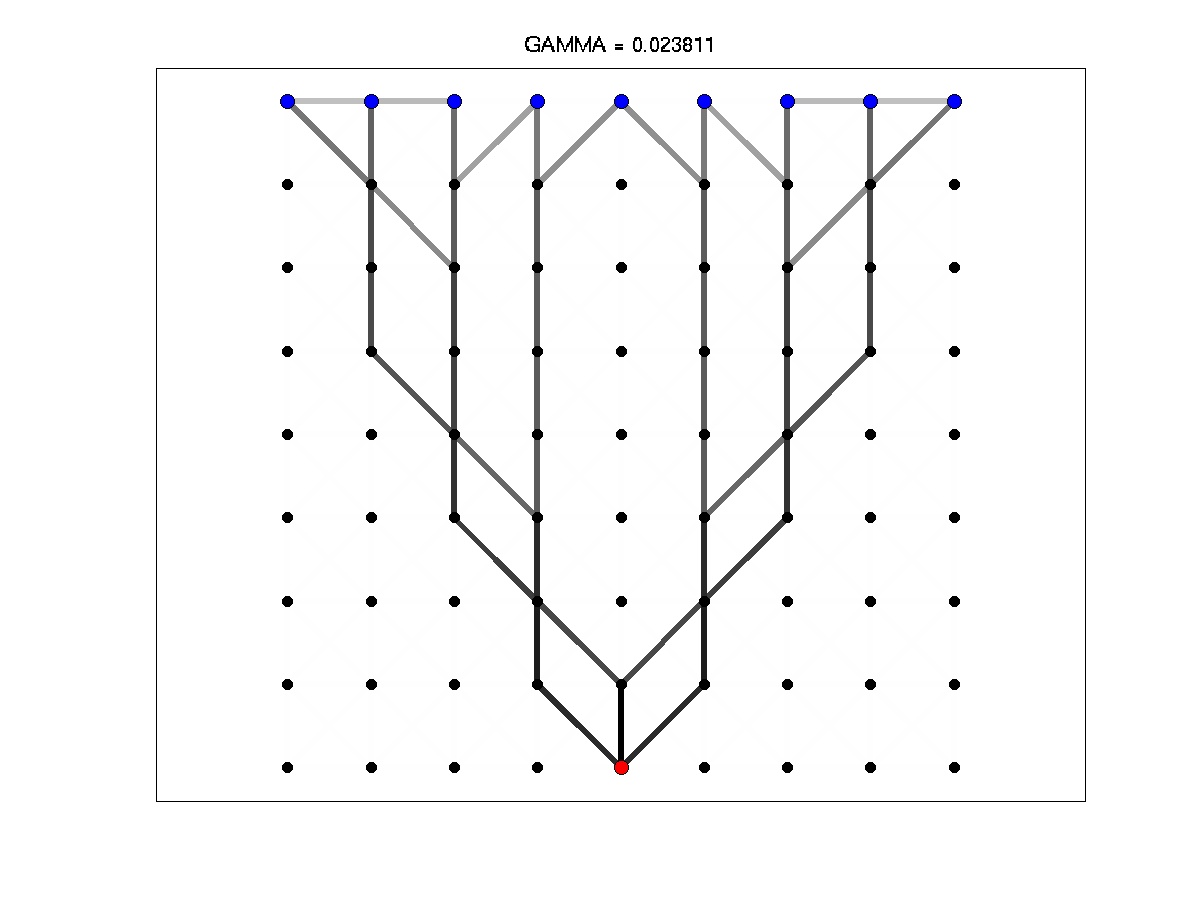}
\hspace{-.6cm} \includegraphics[scale=.113]{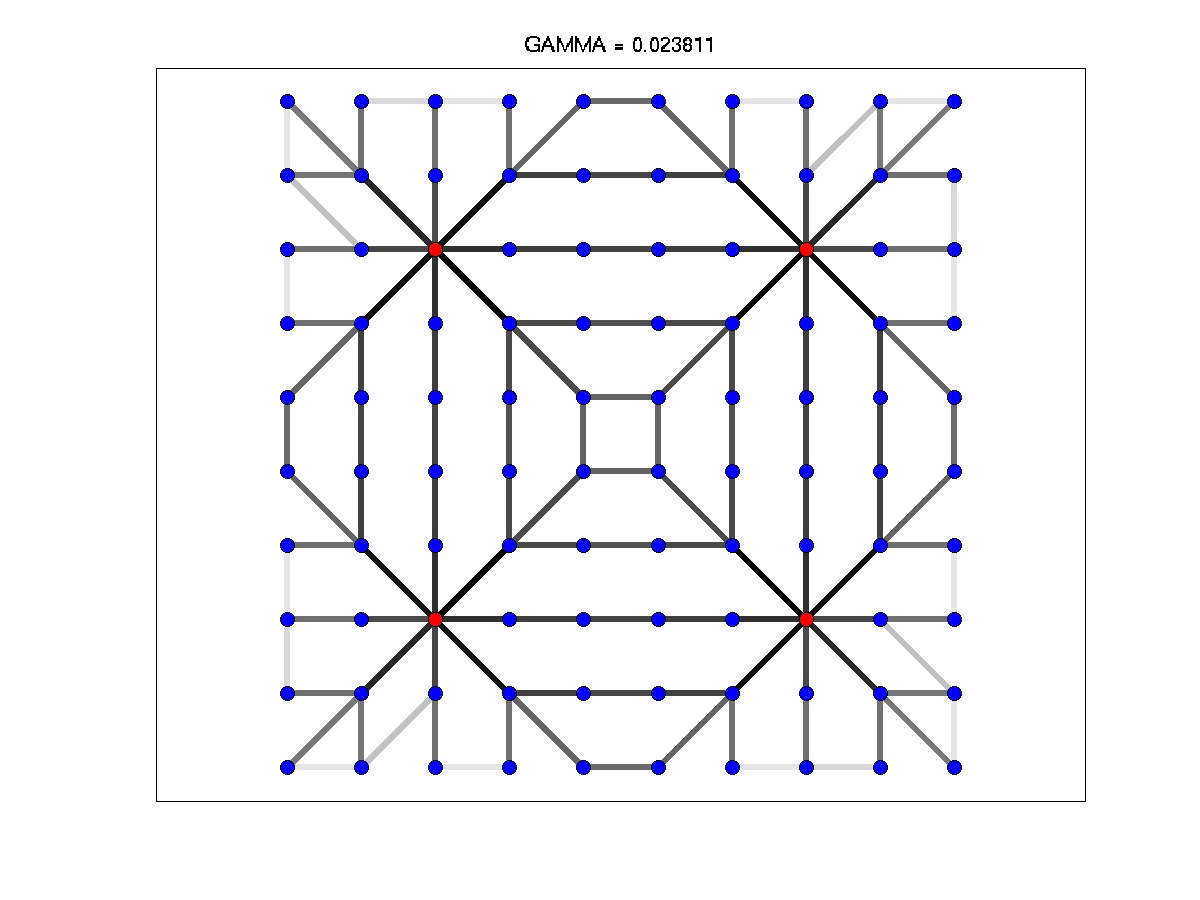}

\vspace{-.3cm}
\caption{Illustration of robust network designs. Compare to Fig.~2.}
\vspace{-.3cm}

\end{figure}

\section{Summary and Future Work}

In summary, we have developed an optimization approach
to design electric power transmission networks with the
aim of balancing network efficiency versus the cost of building
the network.  At the core of our methods lies a convex network
optimization problem generalizing methods of \cite{gbs08,bvg01}.
We also have proposed non-convex extensions of this basic
line-sizing problem to further encourage network sparsity.  This
allows the heuristic design of the network structure $G$ by seeking
a sparse solution within an over-complete graph.   We developed
a heuristic solution technique using smoothed annealing and
majorization-minimization methods \cite{bz87,hl04,cwb08}.  So far,
the experimental results obtained by these methods in toy
problems have yielded very reasonable networks that appear to be
optimal or near-optimal solutions of the proposed optimization
problem.

There are many possible extensions of the basic optimization model
we have developed. We begin by listing some straight-forward extensions
that may give a better fit to real-world applications:
\begin{itemize}

\item Modeling renewable power generation.  This
may be handled as we have treated consumer nodes, but with positive $\bar{b}_i$.

\item Incorporating constraints on the maximum output of generators.
This is especially relevant in the multiple-generator setting.

\item Modeling power storage capabilities (e.g. hybrid and
electric vehicles).  These are nodes of the network that may
absorb power from the network when there is a surplus and then
re-emit this power when demand is high.

\item Allowing for load shedding.  For a number of reasons, it may
become necessary that not all of the demand for power can be met
so that load shedding becomes necessary.  It would be good to
treat this somehow both in our random current model and in how
we model the handling of line and/or generator failures.

\item Putting (convex) constraints on the power dissipation and/or
current per line.  These is important to avoid over-loading lines
in the first place (to avoid cascading failures).

\item Rather than designing networks from scratch, we may also
plan extensions/upgrades to existing networks in a similar manner simply
by including existing line conductances at no cost.

\end{itemize}
A less trivial direction to explore is that of directly
treating the AC power flow problem (rather then using the
leading order DC approximation).  However, so far it is unclear
if this can be usefully treated within a convex optimization framework.

Another direction to explore concerns developing more efficient
algorithms.  The methods we are using so far all involve convex
optimization procedures with per-iteration complexity that grows
essentially as $\mathcal{O}(n^3)$ (fixing the degree of $G$).
We anticipate that more scalable algorithms (e.g., $\mathcal{O}(n^{3/2})$
for planar or near-planar graphs) should be possible using formulations that
introduce auxiliary variables so as to allow Newton's method
to use more efficient linear solvers (e.g. nested dissection) 
that exploit sparsity of the initial graph $G$.

\ignore{Otherwise, it would be interesting to better motivate
the heuristics developed for non-convex optimization from
a theoretical perspective.   Specifically, it would be
useful to identify some sufficient conditions for optimality
of those heuristics (in special cases) or to bound the
sub-optimality of the approximate solutions.  So far, we
have not found a good approach to these questions.}

\section*{Acknowledgments}

We are thankful to all the participants of the ``Optimization and Control for
Smart Grids" LDRD DR project at Los Alamos and Smart Grid Seminar Series at
CNLS/LANL for multiple fruitful discussions, and especially to S.~Backhaus for
critical reading of the manuscript. We also thank R.~Chartrand for a helpful
discussion at the beginning of the project.  Research at LANL was carried out
under the auspices of the National Nuclear Security Administration of the U.S.
Department of Energy at Los Alamos National Laboratory under Contract No. DE
C52-06NA25396. M.~Chertkov acknowledges partial support of NMC via NSF
collaborative grant CCF-0829945 on ``Harnessing Statistical Physics for
Computing and Communications.''

\bibliographystyle{plain}
\bibliography{netopt_cdc}

\end{document}